\documentclass{conm-p-l}

\copyrightinfo{2006}{}

\usepackage{graphicx}

\newtheorem{theorem}{Theorem}[section]
\newtheorem{lemma}[theorem]{Lemma}

\theoremstyle{definition}
\newtheorem{definition}[theorem]{Definition}

\theoremstyle{remark}
\newtheorem{remark}[theorem]{Remark}

\numberwithin{equation}{section}

\newcommand{\e}{\epsilon}

\newcommand{\Ga}{\Gamma}

\newcommand{\dl}{\delta}
\newcommand{\Dl}{\Delta}

\newcommand{\ra}{\rightarrow}
\newcommand{\al}{\alpha}
\newcommand{\be}{\beta}
\newcommand{\sg}{\sigma}
\newcommand{\Sg}{\Sigma}

\newcommand{\pa}{\partial}

\newcommand{\hv}{\hat{v}}

\newcommand{\om}{\omega}

\newcommand{\na}{\nabla}

\newcommand{\tx}{\tilde{x}}

\renewcommand{\O}{{\mathcal O}}

\newcommand{\vphi}{\varphi}

\newcommand{\tv}{\tilde{v}}
\newcommand{\tw}{\tilde{w}}

\begin{document}

\title[The Poincar\'e Recurrence Problem]
{The Poincar\'e Recurrence Problem of Inviscid Incompressible Fluids}

\author{Y. Charles Li}
\address{Department of Mathematics, University of Missouri, 
Columbia, MO 65211, USA}
\curraddr{}
\email{cli@math.missouri.edu}
\thanks{}


\subjclass{Primary 37, 76; Secondary 35, 34}
\date{}

\dedicatory{}


\begin{abstract}
Nadirashvili presented a beautiful example showing that the Poincar\'e recurrence
does not occur near a particular solution to the 2D Euler equation of inviscid 
incompressible fluids. Unfortunately, Nadirashvili's setup of the phase space is not 
appropriate, and details of the proof are missing. This note fixes that.
\end{abstract}

\maketitle










\section{Introduction}

2D Euler equation of inviscid incompressible fluids is the key in understanding 
the chaotic (turbulent) solutions of 2D Navier-Stokes equation in the infinite 
Reynolds number limit \cite{Li04} \cite{LL06}. Moreover, it is the simplest 
fluid equation. There are two distinct phenotypes of chaos: sensitive dependence 
on initial data, and recurrence. Sensitive dependence on initial data can often be 
proved by shadowing technique or Smale horseshoe construction \cite{Li04}. To 
accomplish this, one often needs detailed information on the dynamics (e.g. existence 
of a homoclinic orbit). This poses tremendous analytical challenge for 2D Euler 
equation \cite{LL06}. On the other hand, the well-known Poincar\'e recurrence theorem 
was proved primarily from the first-principle of measure theory under extremely 
general conditions of finite measure space and measure-preserving map. Therefore, it 
seems to have a good chance of success even for 2D Euler equation. The serious 
complication comes from the fact that natural finite dimensional measures (e.g. Gibbs measure) 
do not have good counterparts in infinite dimensions. It is well-known that the 
kinetic energy and enstrophy are invariant under the 2D Euler flow. But it is difficult to 
use them to define finite measures in infinite dimensions. It seems possible to study 
the Poincar\'e recurrence problem directly from Banach norms rather than measures. 
Nadirashvili gave a counter-example of Poincar\'e recurrence along this line 
\cite{Nad91}. Of course, the most exciting and challenging problem shall be the general 
description of Poincar\'e recurrence or non-recurrence for 2D Euler equation directly 
from Banach norms. This is our main research project. In this note, we will fix the problems 
in the Nadirashvili's proof of Poincar\'e non-recurrence near a particular solution to 
the 2D Euler equation. Because Nadirashvili's is the first and a beautiful counter-example 
in this content, we feel that such a note is worth-while.

\section{The Poincar\'e Recurrence Theorem}

\begin{figure}
\includegraphics[width=4.0in,height=3.0in]{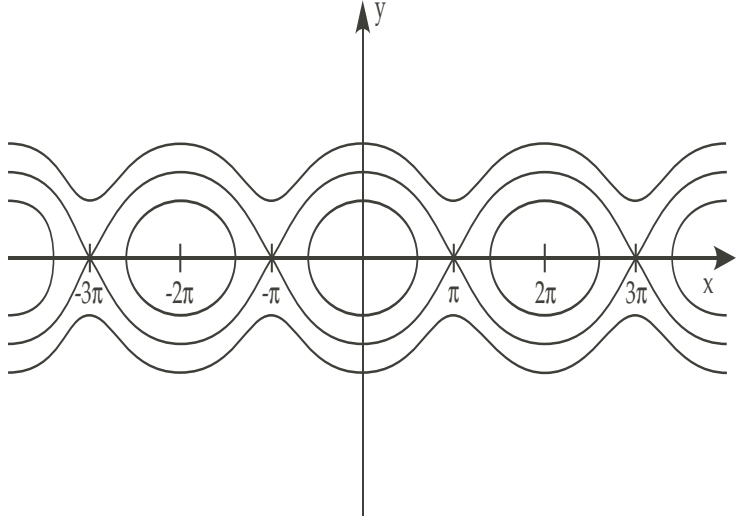}
\caption{The phase plane diagram of the pendulum equation.}
\label{ppd}
\end{figure}

\begin{theorem}
Let ($X, \Sg , \mu$) be a finite measure space and $f\ : \ X \mapsto X$ be a 
measure-preserving transformation. For any $E \in \Sg$ ($\sg$-algebra of 
subsets of $X$), the measure
\[
\mu (\{ x \in E \ | \ \exists N,\  f^n(x) \not\in E \ \forall n > N \} ) = 0 \ .
\]
That is, almost every point returns infinitely often.
\end{theorem}
\begin{proof}
This proof is owned by Koro. 

http://planetmath.org/?op=getobj\&from=objects\&id=6035

Let
\[
A_n = \bigcup_{k=n}^{+\infty}f^{-k}(E) \ ,
\]
then
\[
E \subset A_0\ , \quad A_j \subset A_i \quad \forall i \leq j \ ,
\]
and 
\[
A_j = f^{i-j}(A_i) \ .
\]
Thus
\[
\mu (A_i) = \mu (A_j) \quad \forall i,j \geq 0 \ ,
\]
and
\[
\mu (A_0-A_n) = \mu (A_0) - \mu (A_n) = 0 \ , \quad \forall n \ .
\]
We have
\[
\mu (E-\bigcap_{n=1}^{+\infty}A_n) \leq \mu (A_0-\bigcap_{n=1}^{+\infty}A_n) 
= \mu (\bigcup_{n=1}^{+\infty}(A_0-A_n)) = 0 \ .
\]
Notice that
\[
E-\bigcap_{n=1}^{+\infty}A_n = \{ x \in E \ |\ \exists N,\  f^n(x) \not\in E 
\ \forall n > N \}\ .
\]
The theorem is proved.
\end{proof}
\begin{remark}
The proof in \cite{Wal82} is incorrect. The geometric intuition of the Poincar\'e 
recurrence theorem is that in a finite measure space (or invariant subset), the 
images of a positive measure set under a measure-preserving map will have no room left but 
intersect the original set repeatedly.
\end{remark}
The measure of the space $X$ being finite is crucial. For example, consider 
the two-dimensional Hamiltonian system of the pendulum
\begin{equation}
\dot{x} = y \ , \quad \dot{y} = - \sin x \ . \label{pex}
\end{equation}
Its phase plane diagram is shown in Figure \ref{ppd}. If the invariant region 
includes orbits outside the cat's eyes, then the measure of the region will not 
be finite, and the Poincar\'e recurrence theorem will not hold. One can see clearly 
that the orbits outside the cat's eyes will drift to infinity. The boundaries of the 
cat's eyes are called separatrices (heteroclinic orbits). For 2D Euler equation of 
inviscid incompressible fluids, it has been conjectured that unstable fixed points 
are connected by heteroclinic orbits \cite{LL06}. So the Poincar\'e recurrence 
shall not be sought in the whole phase space.

\section{The Poincar\'e Recurrence Problem of Inviscid Incompressible Fluids}

\subsection{Nadirashvili's Example}

In the case of 2D Euler equation of inviscid incompressible fluids, Nadirashvili was the first 
one who gave an example showing that the Poincar\'e recurrence does not occur. 
Unlike in finite dimensions, here there is no proper measure. Therefore, there 
is no invariant region with a properly defined finite measure. Thus the condition
of the Poincar\'e recurrence theorem does not hold. Even in finite dimensions, the 
finite measure condition is crucial as shown in example \ref{pex}. 

In the original article \cite{Nad91} of Nadirashvili, the setup was incorrect,
where a $C^1$ velocity space was taken as the phase space. The reason is that one 
needs at least $C^{1+\al}$ ($0<\al <1$) initial data in velocity to get $C^1$ 
(in space) velocity solution of the 2D Euler equation \cite{Kat67} \cite{Kat00}. 
In general, $C^{1+\al}$ ($0<\al <1$) initial data can lead to $C^{1+\be}$ 
($0<\be <\al <1$) solutions \cite{Kat00}. So $C^{1+\al}$ ($0<\al <1$) is not a 
good space for a dynamical system study either. On the other hand, Nadirashvili's
is a beautiful example; therefore, it should be set right. In this note, we will select 
an appropriate phase space.

\begin{figure}
\includegraphics[width=4.0in,height=4.0in]{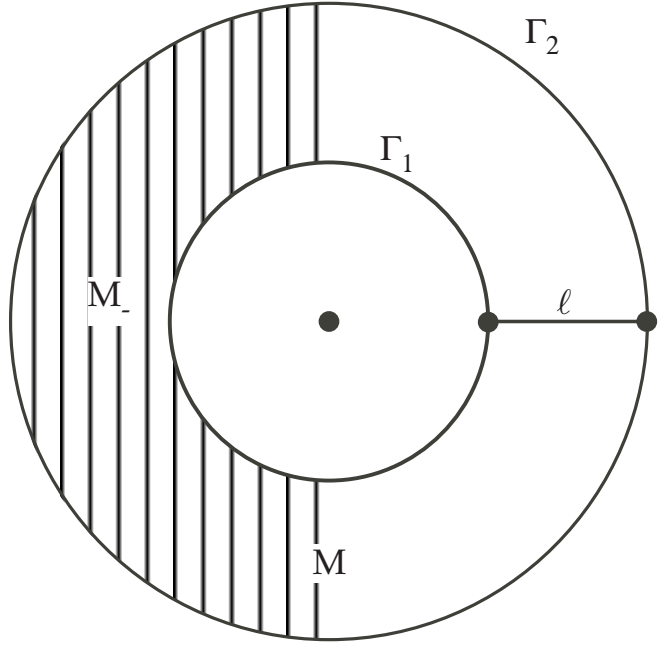}
\caption{The fluid domain for the example of Nadirashvili.}
\label{Nad}
\end{figure}

Let $M$ be the annulus
\[
M = \left \{ x \in \mathbb{R}^2 \ | \ 1 \leq |x| \leq 2 \right \} \ ,
\]
$\Ga_1$ and $\Ga_2$ be the inner and outer boundaries, then $\pa M = \Ga_1 \cup \Ga_2$
(see Figure \ref{Nad}). 
Consider the 2D Euler equation in its vorticity form
\begin{equation}
\left \{ \begin{array}{l} 
\pa_t \om + u \cdot \na \om = 0 \ , \cr 
\om = \text{ curl }u\ , \quad \na \cdot u = 0\ ,\cr 
u \cdot n = 0 \quad \text{on } \pa M = \Ga_1 \cup \Ga_2 \ , 
\quad \int_{\Ga_1} u \cdot dl = \sg_1 \ , \cr
u(0,x)=u^0(x) \ . \cr 
\end{array} \right.
\label{Euler}
\end{equation}
Let $S$ be the space
\[
S = \{ \om \ | \ \om \in C^1(M) \} \ .
\]
Then for any $\om^0\in S$, there is a unique solution to the 2D Euler equation 
(\ref{Euler}) $\om (t) \in S$ for all $t \in \mathbb{R}$, with the initial data 
$\om (0) = \om^0$. In fact, $\om (t)$ is a classical solution in the sense that 
(\ref{Euler}) is satisfied for all $t$ and $x$ (\cite{MP94}, pp.73). In terms of 
the new phase space, the Nadirashvili's theorem \cite{Nad91} can be restated as follows:
\begin{theorem}
There exists a $\xi \in S$, $\e >0$ and $T>0$ such that for any initial condition 
$\om (0) \in S$ satisfying $\| \om (0)-\xi \|_{C^1} < \e $, the corresponding solution 
$\om (t)$ of the 2D Euler equation satisfies $\| \om (t)-\xi \|_{C^1} > \e $ for 
all $t >T$. 
\end{theorem}
We need the following lemma to prove the theorem.
\begin{lemma}
For ($\om , v$) satisfying $\om = \text{ curl }v$, $\na \cdot v=0$, $v\cdot n = 0$
on $\pa M$, $\int_{\Ga_1} v \cdot dl = 0$; one has the inequality
\[
\| v \|_{C^0} \leq C \| \om \|_{C^0} \ .
\]
\label{vvl}
\end{lemma}
A proof of this lemma can be found in (\cite{MP94}, Lemma 3.1, pp.67). Below we 
give a rather detailed proof too.
\begin{proof}
By the Biot-Savart formula for multi-connected domain (\cite{MP94}, pp.16), we have
\[
v=\hv + \tv +\na \vphi
\]
where
\begin{eqnarray*}
& & \hv = \frac{\int_{\Ga_1} v \cdot dl}{2\pi} |x|^{-2} (-x_2, \ x_1) = 0 \ , \\
& & \tv = \frac{1}{2\pi}\int_M |x-\tx |^{-2} \om (\tx ) \left ( -(x_2 -\tx_2), 
\ x_1 -\tx_1 \right )d \tx \ , \\
& & \Dl \vphi = 0 \ , \quad \frac{\pa \vphi}{\pa n} = - \tv \cdot n \ 
\text{ on } \Ga_1 \cup \Ga_2\ , \quad n= \frac{x}{|x|}\ .
\end{eqnarray*}
It is easy to estimate $\| \tv \|_{C^0}$:
\[
\| \tv \|_{C^0} \leq \frac{\| \om \|_{C^0}}{2\pi} \int_M |x-\tx |^{-1}d \tx 
\leq C \| \om \|_{C^0} \ .
\]
Estimating $\| \na \vphi \|_{C^0}$ is more complicated. We know that $\vphi$ is 
given by the single layer potential (\cite{Fol76}, pp.171)
\begin{equation}
\vphi (x) = \int_{\Ga_1 \cup \Ga_2} N(x, \tx ) f(\tx ) d \tx \ , 
\quad N(x, \tx )=\frac{1}{2\pi} \ln |x-\tx | \ ,
\label{slp}
\end{equation}
where the moment $f$ is given by (\cite{Fol76}, pp.174)
\begin{equation}
-\frac{1}{2} f(x) + \int_{\Ga_1 \cup \Ga_2} \pa_{n_x}N(x, \tx )f(\tx ) d \tx 
=-\tv \cdot n \ .
\label{moe}
\end{equation}
Using an approximation (e.g. as in \cite{Fol76}, pp.95), one can show that 
\[
\na \cdot \tv = 0 \quad \text{for all } x \in \mathbb{R}^2 \ .
\]
Thus
\[
\int_{\Ga_j} n \cdot \tv \ dx = 0 \quad \text{for } j=1,2,
\]
which is an if and only if condition for the existence of a solution to (\ref{moe}) 
in $L^2(\Ga_1 \cup \Ga_2)$, and
\[
\| f\|_{L^2(\Ga_1 \cup \Ga_2)} \leq C \| \tv \cdot n \|_{L^2(\Ga_1 \cup \Ga_2)}
\leq C \| \om \|_{C^0} \ .
\]
In fact $f \in C(\Ga_1 \cup \Ga_2)$ (\cite{Fol76}, pp.160 (3.14)). Notice that
(\cite{Fol76}, pp.163)
\[
\pa_{n_x}N(x, \tx ) = \frac{1}{2\pi} \frac{(x-\tx )\cdot n_x}{|x-\tx |^2} 
= \frac{1}{2\pi} \frac{\cos (x-\tx , n_x )}{|x-\tx |} \ ,
\]
and
\[
\cos (x-\tx , n_x ) = \O (|x-\tx |) \quad \text{as } \tx \ra x \text{ on }
\Ga_1 \cup \Ga_2 \ .
\]
Thus $\pa_{n_x}N(x, \tx )$ is bounded on $(\Ga_1 \cup \Ga_2) \times (\Ga_1 \cup \Ga_2)$.
We have from (\ref{moe}) that
\[
\| f\|_{C^0(\Ga_1 \cup \Ga_2)} \leq C \| f\|_{L^2(\Ga_1 \cup \Ga_2)} +
\| \tv \cdot n \|_{C^0(\Ga_1 \cup \Ga_2)} \leq C \| \om \|_{C^0} \ .
\]
This is the estimate we need. From (\ref{slp}),
\[
\na \vphi (x) = \frac{1}{2\pi} \int_{\Ga_1 \cup \Ga_2} \frac{x-\tx }{|x-\tx |^2} 
f(\tx ) d \tx \ .
\]
Thus
\[
\| \na \vphi \|_{C^0} \leq C \| f\|_{C^0(\Ga_1 \cup \Ga_2)} \leq C \| \om \|_{C^0} \ .
\]
The proof is complete. 
\end{proof}
\begin{remark}
This lemma was not proved in the original article \cite{Nad91}. In \cite{AK98}, the authors 
commented on the possibility of a proof by a maximal principle, but no detail was available.
\end{remark}
\begin{proof} [Proof of the Theorem]
Let 
\[
u^* = |x|^{-2}(-x_2, \ x_1)\ .
\]
Then
\begin{eqnarray*}
& & \text{curl }u^* = \na \cdot u^* = 0\ ,\\ 
& & u^* \cdot n = 0 \quad \text{on } \pa M = \Ga_1 \cup \Ga_2 \ , \\
& & \int_{\Ga_1} u^* \cdot dl = 2\pi \ .
\end{eqnarray*}
Let 
\[
M_- = \{ x \in M \ | \ x_1 <0 \} \ , \quad \ell =  \{ x \in M \ | \ x_2=0, x_1 >0 \} \ .
\]
Choose $\xi$ such that
\[
\| \xi \|_{C^1} < 4\e \ , \quad \xi |_{M_-} = 0 \ , \quad \xi |_{\ell} > 2 \e \ .
\]
For any initial condition $\om (0) \in S$ such that $\| \om (0)-\xi \|_{C^1} < \e$,
let $\om (t)$ be the corresponding solution of the 2D Euler equation (\ref{Euler})
with $\int_{\Ga_1}u \cdot dl = \int_{\Ga_1} u^* \cdot dl$. Let $v$ be the 
corresponding velocity given by Lemma \ref{vvl}. When $\e$ is small enough, 
$\| v\|_{C^0} < \frac{1}{4}$ for all $t \in \mathbb{R}$. Let $u = u^* +v$. Then 
($\om , u$) solves the 2D Euler equation (\ref{Euler})
with $u(0,x) = u^*(x) +v(0,x)$, $\int_{\Ga_1}u \cdot dl = \int_{\Ga_1} u^* \cdot dl$.
$u(t)$ defines for each $t\in \mathbb{R}$ an area-preserving diffeomorphism $g^t$ of 
the annulus $M$ (\cite{MP94}, Theorem 3.1, pp.72). The boundaries $\Ga_1$ and $\Ga_2$
are mapped into themselves by $g^t$. We have
\[
u(t)|_{\Ga_1} > 3/4\ , \quad u(t)|_{\Ga_2} < 3/4 \quad \text{in polar coordinate}.
\]
Thus the corresponding angular velocities are greater than $3/4$ on $\Ga_1$ and 
smaller than $3/8$ on $\Ga_2$. The image $\ell_t = g^t(\ell )$ of the segment 
$\ell$ under the diffeomorphism $g^t$ will wrap around in the annulus with 
angular coordinates of the two ends diverging faster than $\frac{3t}{8}$. So 
$\ell_t$ will wrap around the inner circle with more and more loops. Thus when $t > 
\frac{8\pi}{3}$, $\ell_t$ will always intersect $M_-$. Notice that the value of 
$\om (t)$ is carried 
over by $\ell_t$ and $\om (t) > \e$ on $\ell_t$. Hence, for $t > 
\frac{8\pi}{3}$, we have $\| \om (t) - \xi \|_{C^1} > \e$. 
\end{proof}
\begin{remark}
The setup of $\text{curl } u |_\ell > \dl /4$ v.s. $\e = \dl /4$ in (\cite{AK98}, 
pp.98) is not compatible for the later argument --- a trivial problem that can be 
easily fixed as above. 
\end{remark}

\subsection{Shnirelman's Theorem}

In terms of Lagrangian coordinates, unlike the Nadirashvili's example; 
Shnirelman proved a theorem on the wandering nature of the configuration map 
induced by the 2D inviscid incompressible fluid motion in $\mathbb{T}^2$.
\cite{Shn97}. The form of the fluid (particle trajectory) equation considered by Shnirelman is 
\begin{equation}
\frac{d}{dt} (g, \om ) = \left ( \text{rot}^{-1} (\om \circ g^{-1})\circ g , 0 \right )\ ,
\label{flue}
\end{equation}
defined on $\mathbb{T}^2$, where $g\ : \ \mathbb{T}^2 \mapsto \mathbb{T}^2$ is the configuration map 
induced by the fluid motion in $\mathbb{T}^2$, and $\om$ is the vorticity. 

\begin{definition}
Let $w(x) \in L^2(\mathbb{T}^2)$ and $\tw (\xi )$ ($\xi \in \mathbb{Z}^2$) be its 
Fourier transform. Besov space $B_s$ is the space of functions $w(x)$ with finite norm
\[
\| w\|_{B_s}^2 =\sup_{k \geq 0} \left \{ 2^{2ks} \sum_{2^k \leq |\xi | < 2^{k+1}}
|\tw (\xi )|^2 \right \}\ .
\]
By $(g, \om ) \in X_s = DB_s \times B_{s-1}$, it means that $g(x)-x \in B_s$ and 
$\om (x) \in B_{s-1}$. Denote by $G_t$ the evolution operator of the 2D Euler equation,
i.e. $G_t (g(0), \om (0)) = (g(t), \om (t))$.
\end{definition}
\begin{theorem}
For $s >0$, there exists an open and dense set $Y_s \subset X_s$, such that for each point
$(g, \om ) \in Y_s$ there is a wandering neighborhood $U$ of $(g, \om )$, i.e. there 
exists a $T>0$ (depending on $(g, \om )$ and $U$) such that $G_t (U) \cap U = \emptyset$
for all $t$ such that $|t|>T$. 
\end{theorem}
\begin{remark}
The Besov space $B_s$ differs only slightly from the Sobolev spaces: 
$H^s \subset B_s \subset H^{s-\e}$ for any $s$ and $\e >0$. The smooth functions are 
not dense in the Besov space $B_s$ (unlike the Sobolev space $H^s$). For $s>3$, the 
fluid equation (\ref{flue}) is globally wellposed in $X_s$. The configuration map $g$ can be viewed as 
the family of all fluid trajectories. The theorem says that most of the families of 
all fluid trajectories are wandering. A more challenging problem will be the general 
description of wandering solutions of the 2D Euler equation on $\mathbb{T}^2$.
\end{remark}


\begin{thebibliography}{99}

\bibitem{AK98} V. Arnold, B. Khesin, {\it Topological Methods in Hydrodynamics}, 
Springer, Applied Math. Sci., vol.125, (1998), pp.98.

\bibitem{Fol76} G. Folland, {\it Introduction to Partial Differential Equations}, 
Princeton University Press, (1976).

\bibitem{Kat67} T. Kato, On classical solutions of the two-dimensional 
non-stationary Euler equation, {\it Arch. Rat. Mech. Anal.} {\bf 25} (1967), 188-200.

\bibitem{Kat00} T. Kato, On the smoothness of trajectories in incompressible 
perfect fluids, {\it Contemp. Math.} {\bf 263} (2000), 109-130.

\bibitem{LL06} Y. Lan, Y. Li,  On the dynamics of Navier-Stokes and Euler equations, 
{\it Submitted} (2006).

\bibitem{Li04} Y. Li, {\it Chaos in Partial Differential Equations}, International Press,
(2004).

\bibitem{MP94} C. Marchioro, M. Pulvirenti, {\it Mathematical Theory of Incompressible 
Nonviscous Fluids}, Springer, Applied Math. Sci., vol.96, (1994).

\bibitem{Nad91} N. Nadirashvili, Wandering solutions of Euler's D-2 equation,
{\it Funct. Anal. Appl.} {\bf 25, no.3} (1991), 220-221.

\bibitem{Shn97} A. Shnirelman, Evolution of singularities, generalized Liapunov 
function and generalized integral for an ideal incompressible fluid,
{\it Amer. J. Math.} {\bf 119} (1997), 579-608.

\bibitem{Wal82} P. Walters, {\it An Introduction to Ergodic Theory}, Springer-Verlag, 
(1982), pp.26.



\end{thebibliography}
\end{document}